\documentclass[1p]{elsarticle}

\usepackage{hyperref}

\usepackage{amssymb}
\usepackage{amstext}
\usepackage{amsthm}
\usepackage{mathtools}
\usepackage{enumerate}  

\newtheorem{definition}{Definition}
\newtheorem{theorem}{Theorem}
\newtheorem{prop}{Proposition}

\newtheorem{corollary}{Corollary}

\newtheorem{question}{Question}

\newcommand{\cA}{\mathcal{A}}
\newcommand{\cG}{\mathcal{G}}
\newcommand{\cH}{\mathcal{H}}

\DeclareMathOperator{\tr}{tr}

\journal{European Journal of Combinatorics}

\allowdisplaybreaks

\begin{document}
\begin{frontmatter}

\title{Characteristic Power Series of Graph Limits}

\author[mymainaddress]{Joshua N.~Cooper}
\address{Department of Mathematics, University of South Carolina, Columbia, SC USA}
\cortext[mycorrespondingauthor]{Corresponding author}
\ead{cooper@math.sc.edu}

\begin{abstract}
In this note, we show how to obtain a ``characteristic power series'' of graphons -- infinite limits of dense graphs -- as the limit of normalized reciprocal characteristic polynomials.  This leads to a new characterization of graph quasi-randomness and another perspective on spectral theory for graphons, a complete description of the function in terms of the spectrum of the graphon as a self-adjoint kernel operator.  Interestingly, while we apply a standard regularization to classical determinants, it is unclear how necessary this is.
\end{abstract}


\end{frontmatter}

\section{Introduction}

A research direction began in the 1980's with graph quasi-randomness, extended through the 1990's and early 2000's with generalizations to non-uniform graph distributions and other combinatorial objects, became graph limit theory in the mid-2000's, and culminated in Lov\'asz's now-canonical text \cite{Lo12}.  The central idea is that, if a sequence of graphs $G_n$ with number of vertices tending to infinity has the property that the density of any particular subgraph tends to a limit, then $G_n$ itself tends to a limit object $\cG$, called a ``graphon''.  There are several mutually (though non-obviously) equivalent ways to view graphons, and a central one is as a self-adjoint kernel operator from $L^1([0,1])$ to $L^\infty([0,1])$, an object type for which a well-established spectral theory exists.  In particular, a graphon, when thought of as this kernel, is a symmetric function $[0,1] \times [0,1] \rightarrow [0,1]$, with two such kernels $U$ and $W$ ``weakly isomorphic'' (i.e., giving rise to the same graphon) if there are measure-preserving maps $f,g : [0,1] \rightarrow [0,1]$ so that $U(f,f) = W(g,g)$ almost everywhere. Indeed, we use the same notation throughout for $\cG$ as well as (any representative of the equivalence class of) its kernel. There is also a key notion of subgraph density for graphons, with the property that densities in convergent sequences of graphs in a sequence converge to their densities in the limit graphon, often referred to as ``left-convergence''.  Furthermore, Szegedy (\cite{Sze11}) introduced a spectral theory of graphons by studying the eigenpairs of their kernels and showed that it is a natural analogue of the adjacency spectral theory of finite graphs.  Here, we extend this perspective by showing that graphons are associated with a power series which is a certain normalized limit of the characteristic polynomial of graphs.  Furthermore, we give an equivalent definition of this ``characteristic power series'' $\psi_\cG(z)$ of a graphon $\cG$ via a regularized determinant of its corresponding kernel function.

We also show that the characteristic power series can be used to characterize ``quasi-randomness''.  Suppose that $\{G_n\}_{n \geq 1}$ is a sequence of graphs with $|V(G_n)| = n$.  (In truth, all that is needed is that $|V(G_n)| \rightarrow \infty$, but this is not more general.)  We write $G = G_n$ for simplicity and $(\#H \subseteq G)$ for the number of labelled, not-necessarily induced copies of $H$ as subgraphs in $G$ (i.e., injective homomorphisms from $H$ to $G$).  Then, by a classic 1989 paper of Chung, Graham, and Wilson (\cite{ChuGraWil89}), there is a large set of random-like properties (properties which hold asymptotically almost surely for graphs in the Erd\H{o}s-R\'enyi model $G(n,p)$) which are mutually equivalent, and are therefore collectively referred to as (the sequence of graphs) $G$ being ``quasi-random''\footnote{We state the version of this theorem for all densities $p \in [0,1]$, but it actually appeared first only for $p=1/2$.}. Namely, let 
\begin{itemize}
    \item $P_1(s)$ denote the property that the number of labelled occurrences of each graph on $s$ vertices as an induced subgraph of $G$ is $n^s (p^{|E(H)|} (1-p)^{\binom{s}{2}-|E(H)|} + o(1))$
    \item $P'_1(s)$ denote the property that $(\# H \subset G) = n^s (p^{|E(H)|} + o(1))$ for each graph $H$ on $s$ vertices
    \item $P_2(t)$ denote the property that $|E(G)| \geq (p+o(1)) n^2/2 $ and $(\# C_t \subseteq G) \leq (1+o(1)) (np)^t$, where $C_t$ is the $t$-cycle
    \item $P_3$ denote the property that $|E(G)| \geq (p+o(1))n^2/2$, $\lambda_1 = pn (1+o(1))$, and $\lambda_2 = o(n)$, where $|\lambda_1| \geq \cdots \geq |\lambda_n|$ are the complete set of adjacency eigenvalues of $G$ with multiplicity
    \item $P_4$ denote the property that, for all $S \subseteq V(G)$, $E(G[S]) = p |S|^2/2 + o(n^2)$, where $G[S]$ denotes the subgraph of $G$ induced by $S$
    \item $P_5$ denote the property that, for all $S \subseteq V(G)$ with $|S| = \lfloor n/2 \rfloor$, $E(G[S]) = p n^2/8 + o(n^2)$, where $G[S]$ denotes the subgraph of $G$ induced by $S$
    \item $P_6$ denote the property that
    $$
    \sum_{v,w \in V(G)} \left | |N_G(v) \cap N_G(w)| - np^2 \right | = o(n^3)
    $$    
\end{itemize}

\begin{theorem}[Chung-Graham-Wilson \cite{ChuGraWil89}] \label{thm:QRG} For $s \geq 4$ arbitrary and $t \geq 4$ even, and any fixed $p \in [0,1]$,
$$
P_2(4) \Leftrightarrow P_2(t) \Leftrightarrow P_1(s) \Leftrightarrow P'_1(s) \Leftrightarrow P_3 \Leftrightarrow P_4 \Leftrightarrow P_5 \Leftrightarrow P_6
$$
\end{theorem}

If a graph sequence $G$ has these properties, it is called {\em $p$-quasi-random}, and these properties and any others also equivalent to them are known as {\em ($p$-)quasi-random properties}.  Many other quasi-random properties have been added since to the list above, such as other families $\mathcal{F}$ of graphs whose occurrence as subgraphs at the ``random-like rate'' implies these properties (note that $\{K_2,C_4\}$ is the ``forcing'' family given by $P_2(4)$), and also that $G_n$ converges to a $p$-constant graphon.

Here, we propose to add another property.  First, given a left-convergent sequence of graphs $G$, let $\cG$ be their limit, and let $\phi_n \in \mathbb{C}[x]$ denote the (adjacency) characteristic polynomials of $G_n$.  Recall that the characteristic polynomial of a graph $G$ with adjacency matrix $A$ is defined to be $\phi(x) = \det(Ix - A)$; it is easy to see that $\phi$ is monic and has only real roots.  Define
$$
\psi_\cG(z) = \left \{ \begin{array}{ll} \lim_{n \rightarrow \infty} \left ( x^{-n} \phi_{n}(x) \right ) \Big |_{x = n/z} & \textrm{ if } z \neq 0\\
1 & \textrm{ otherwise },
\end{array} \right .
$$
if the (pointwise) limit exists.  We call $\psi_\cG(z)$ the ``characteristic power series'' of the graphon $\cG$, and it is essentially a normalized limit of the reciprocal polynomials of the characteristic functions of $G$.  In Theorem \ref{thm:convergence} below, we show that $\psi_\cG$ is indeed well-defined, i.e., independent of the sequence of graphs left-converging to $\cG$.

\section{Characteristic Polynomials of Large Graphs}

The following classical result will be useful in describing the coefficients of $\psi_\cG(z)$.  Let $\mathcal{H}_k$ -- sometimes called ``elementary graphs'' -- denote the family of unlabelled graphs on $k$ vertices each of whose components is an edge or a cycle.

\begin{theorem}[Harary-Sachs \cite{Ha62}] \label{thm:HS}  Suppose $G$ is a graph on $n$ vertices, and $k \geq 0$ is an integer.  The coefficient of $x^{n-k}$ in $\phi_G(x)$ is
$$
\sum_{H \in \cH_k} (-1)^{c(H)} 2^{z(H)} [\# H \subseteq G] 
$$
where $c(H)$ is the number of components of $H$, $z(H)$ is the number of cycles of $H$, and $[\# H \subseteq G]$ denotes the number of subsets of edges of $G$ which are isomorphic to $H$.  (When $k=0$, $\cH_k$ is the singleton consisting only of the empty graph $\epsilon$; take $c(\epsilon) = z(\epsilon) = 0$.)
\end{theorem}

Clearly, if $H \in \cH_k$, then $|E(H)| = k - c(H) + z(H)$.  Also,
$$
(\# H \subseteq G) = [\# H \subseteq G] \cdot 2^{z(H)} \prod_{i} a_i^{m_i} m_i!
$$
where $H$ has $m_i$ components of size $a_i$ for each $i$.  For simplicity, for a partition $\lambda$ where part $b_i$ occurs $m_i$ times for each $i$ (the $b_i$ all distinct), the quantity $\eta(\lambda)$ is defined by
$$
\eta(\lambda) := \prod_i b_i^{m_i} m_i!
$$
and we write $\eta(H) = \eta(\lambda(H))$ where $\lambda(H)$ is the integer partition of $|E(H)|$ given by the component cardinalities of $H$.  If $\lambda$ is a partition of $n$, then the number of partitions of an $n$-set with structure $\lambda$ (a partition with $m_i$ parts of distinct sizes $b_i$) is given by 
$$
\frac{n!}{b_1!^{m_1} \cdots b_t!^{m_t} m_1! \cdots m_t!} = \frac{n!}{\eta(\lambda) \prod_{i=1}^t (b_i-1)!^{m_i}}
$$
Denote by $\Lambda'_{n,k}$ the set of partitions of $n$ into $k$ parts, each of which is of size at least $2$; for $\lambda \in \Lambda'_{n,k}$, denote its $i$-th largest part by $\lambda_i$, the $i$-th largest integer which appears as a part by $b_i$, and the multiplicity of $b_i$ by $m_i$.

We wish to transform the characteristic polynomial of a graph into a power series so that, as $n \rightarrow \infty$, the limit exists if $G_n$ converges to a graphon.  This entails taking the reciprocal of $\phi_{G_n}(x)$, normalizing the coefficients with appropriate powers of $n$, and then letting $n$ tend to infinity.  So, setting $z = n/x$, we have by Theorem \ref{thm:HS},
\begin{align*}
x^{-n} \phi_{G_n}(x) &= x^{-n} \sum_{k=0}^n x^{n-k} \sum_{H \in \cH_k} (-1)^{c(H)} 2^{z(H)} [\# H \subseteq G_n] \\
&= \sum_{k=0}^n \sum_{\lambda \in \Lambda'_{n,k}} x^{-k} (-1)^{k} \frac{(\# \bigsqcup_i C_{\lambda_i} \subseteq G_n)}{\eta(\lambda)} \\
&= \sum_{k=0}^n \sum_{\lambda \in \Lambda'_{n,k}} z^{k} (-1)^{k} \frac{(\# \bigsqcup_i C_{\lambda_i} \subseteq G_n)}{n^k  \prod_i b_i^{m_i} m_i!},
\end{align*}
where $C_2$ denotes a single edge.  Note that (treating this as a polynomial to avoid defining $0^0$), when $z=0$, the above expression equals $1$ because $\phi_{G}(x)$ is monic.  Denote this polynomial by $\psi_n(z)$.  Write $t(H,G)$ for the ``homomorphism density'' of the $k$-vertex graph $H$ in the $n$-vertex $G$, i.e., the number of (not necessarily injective) homomorphisms from $H$ to $G$ over $n^k$.  Writing $\mu_i$ for the $i$-th eigenvalue of $A(G_n)/n$, we may bound
\begin{align*}
| \psi_n(z) | &\leq \sum_{k=0}^n |z|^k \sum_{\lambda \in \Lambda'_{n,k}} \frac{t(\bigcup_i C_{\lambda_i}, G_n)}{\prod_i b_i^{m_i} m_i!} \\
&= \sum_{k=0}^n \sum_{\lambda \in \Lambda'_{n,k}} \prod_{i=1}^k \frac{|z|^{m_i} t(C_{b_i}, G_n)^{m_i}}{b_i^{m_i} m_i!} \\
&= \prod_{b \geq 2} \sum_{m \geq 0} \frac{|z|^{m} t(C_b, G_n)^{m}}{b^{m} m!} \\
&= \prod_{b \geq 2} \exp \left ( \frac{|z| t(C_{b}, G_n)}{b} \right ) =  \exp \left ( \sum_{b \geq 2} \frac{|z| t(C_{b}, G_n)}{b} \right ).
\end{align*}
Then, since $t(C_{b}, G_n) = \sum_i \mu_i^b$ for each $b \geq 2$,
\begin{align*}
| \psi_n(z) | &\leq \exp \left ( \sum_{b \geq 2} \frac{|z| \sum_i \mu_i^b}{b} \right ) = \exp \left ( - |z| \sum_i [\log (1-\mu_i) + \mu_i ] \right ) \\
&= \left [ \prod_i \left ( 1-\mu_i \right ) \right ]^{-|z|} \cdot \exp \left ( - |z| \sum_i \mu_i \right ) = \left [ \prod_i \left ( 1-\mu_i \right ) \right ]^{-|z|}
\end{align*}
since $\sum_i \mu_i = \tr A(G_n)/n = 0$.  The above product converges if $\sum_i |\mu_i|$ does, which is $n^{-1}$ times the so-called ``energy'' $\mathcal{E}(G_n)$ of $G_n$, the sum of its adjacency singular values.  Since $\mathcal{E}(G_n)$ can be as large as $C n^{3/2}$, we should not hope for $\phi_n(z)$ always to converge.  Indeed, $n^{3/2}$ tends to be the order of magnitude of the energy of dense graphs, i.e., graphs with $\Omega(n^2)$ edges, the only graphs converging to a nontrivial graphon; see, for example, \cite{Ni07}.  However, this is not {\em always} the case: indeed, $\mathcal{E}(K_n) = 2n-2$.

\section{Characteristic Power Series of Graphons}

In order to circumvent convergence problems discussed in the previous section, we introduce the so-called ``regularized characteristic determinant'' $\det^{(p)}(\cA)$ of a linear operator $\cA$:

\begin{definition}
For a positive integer $p$, the {\em regularized characteristic determinant} of a linear operator $\cA$ with discrete spectrum is defined as
$$
\det^{(p)}(I - z \cA) = \prod_j \left [(1 - \lambda_j z) \exp \left ( \sum_{k=1}^{p-1} \lambda_j^k z^k/k \right ) \right ]
$$
where $\lambda_j$ varies over the eigenvalues of $\cA$.
\end{definition}

We then use this definition -- albeit only the $p=2$ case, a.k.a.~the Hilbert-Carleman determinant, after first considering $p=4$ -- for reasons which will be apparent below, to define a characteristic power series of graphons:

\begin{definition} \label{def:charpowerseries}
The {\em characteristic power series} of a graphon $\cG$ is defined by
$$
\psi_\cG(z) =  \det^{(2)}(I - z\cG)  \exp \left ( z^2 \cdot \frac{\|\cG\|_2^2 - \|\cG\|_1}{2} \right ).
$$
\end{definition}

By \cite{GoKr69} (Chapter IV, Section 2), the function $\det^{(p)}(I - z \cA)$ is well-defined and entire (of genus $p-1$) for operators $\cA$ in $\mathfrak{S}_p$, the operators which are Schatten $p$-class, i.e., for which the Schatten $p$-norm $(\sum_i \sigma_i^p)^{1/p}$ is finite, where $\sigma_i$ are the singular values of $\cA$, defined to be the eigenvalues of $\sqrt{\cA^\ast \cA}$.  Since graphons give rise to self-adjoint operators, we will have throughout that $\sigma_i = |\lambda_i|$.  Note that Schatten $2$-class bounded operators are the same as Hilbert-Schmidt operators, which all graphons' corresponding integral transforms are; and Schatten $1$-class are the nuclear or trace-class operators, in which case $\det^{(1)}(I - z \cA)$ is the classical Fredholm determinant and $\tr(\cA) = \sum_j \lambda_j$ is the (signed) trace.  It also follows from \cite{GoKr69} (see Theorem IV.2.1) that $\det^{(p)}(I - z \cA)$ is continuous (uniform convergence on compact sets) with respect to convergence in $p$-norm of $\cA$.

We now present our main theorem, demonstrating that $\psi_\cG$ is indeed well-defined and is an entire function of Laguerre-P\'olya class, i.e., a holomorphic function which is locally the limit of a series of polynomials whose roots are all real.  Laguerre-P\'olya functions have played a prominent role in the study of distributions of zeros of real polynomials and real entire functions (e.g., \cite{BaCrCs01}), early 20th-century attempts to prove the Riemann hypothesis and a recent revival of such methods (see \cite{GrOnRoZa19}), and classical complex analysis.  The fact that $\psi_\cG(z)$ is Laguerre-P\'olya class implies that it has a Hadamard product expression (see, e.g., \cite{HiWi55} p.~42--47):
\begin{equation} \label{eq:hadamard}
\psi_\cG(z) = z^m \exp(a + bz + cz^2) \prod_{r} \left (1-\frac{z}{r} \right )\exp \left (\frac{z}{r} \right)
\end{equation}
where $m$ is a nonnegative integer; $b$ and $c$ are real with $c \leq 0$; and $r$ ranges over the nonzero zeros of $\psi_\cG(z)$.  Note that the definition of $\psi_\cG(z)$ is almost in this form already.  In particular, $m=a=b=0$, and $c = (\|\cG\|_2^2 - \|\cG\|_1)/2$, and the product ranges over the reciprocals $r$ of the nonzero eigenvalues of $\cG$.  

\begin{theorem} \label{thm:convergence} Suppose the graphs $G_n$ converge to the graphon $\cG$.  For the sequence of functions $\psi_n(z)$ corresponding to the sequence of graphs $G_n$:
\begin{enumerate}
    \item $\psi_n(z)$ converges pointwise as $n \rightarrow \infty$.
    \item $\psi_n(z)$ converges uniformly on compact sets as $n \rightarrow \infty$.
    \item Each coefficient of $\psi_n(z)$ converges as $n \rightarrow \infty$.
\end{enumerate}
Furthermore, the limit is $\psi_\cG(z)$, is entire of Laguerre-P\'olya class, and its roots are the reciprocals of the nonzero eigenvalues of $\cG$ (as a self-adjoint kernel operator) with multiplicity.
\end{theorem}
\begin{proof} 
Let $A'_n = A(G_n)/n$, so that
\begin{align*}
    \psi_n(z) &= \left .\left ( x^{-n} \det(xI-A(G_n)) \right ) \right |_{x = n/z} \\
    &= \left .\left ( \det(I-A(G_n)x^{-1}) \right ) \right |_{x = n/z} \\
    &=  \det(I-A'_nz).
\end{align*}
There is a natural definition of the homomorphism densities $t(F,\cG)$ in terms of a certain integral which is standard for graphons, and \cite{Lo12} (7.22) shows that $\sum_j \lambda_j(\cG)^k = t(C_k,\cG)$ for $k \geq 3$ and, by left-convergence, $t(C_k,G_n) \rightarrow t(C_k,\cG)$ for $k \geq 3$.  Thus, in particular, $\sum_j \lambda_j(\mathcal{G})^4$ converges, so $\mathcal{G}$ is Schatten $4$-class and $\det^{(4)}(I-z\mathcal{G})$ is entire. That $G_n \rightarrow \mathcal{G}$ implies that there exists a sequence of measure-preserving maps $\varphi_n$ of $[0,1]$ to itself so that, if $\mathcal{A} : [0,1]^2 \rightarrow [0,1]$ is the kernel defined by $\mathcal{A}_n(x,y) = A(G_n)_{\lceil xn \rceil,\lceil yn \rceil}$, then $\tilde{\mathcal{A}}_n(x,y) := \mathcal{A}_n(\varphi(x),\varphi(y))$ converges to $\mathcal{G}$ in cut norm; note that $\tilde{\mathcal{A}}_n$ has the same spectrum as $A'_n$.  By \cite{Lo12} Lemma 8.12, that $\tilde{\mathcal{A}}_n \rightarrow \mathcal{G}$ in cut norm implies that $\tilde{\mathcal{A}}_n \rightarrow \mathcal{G}$ in Schatten $4$-norm.  As mentioned following Definition \ref{def:charpowerseries}, Theorem IV.2.1 of \cite{GoKr69} then implies that 
\begin{equation} \label{eq4}
\prod_j \left (1 - \lambda_j(A'_n) z \right ) \exp \left ( \sum_{k=1}^{3} \frac{\lambda_j(A'_n)^k z^k}{k} \right ) \! \rightarrow \prod_j \left (1 -z \lambda_j(\cG) \right ) \exp \left ( \sum_{k=1}^{3} \frac{\lambda_j(\cG)^k z^k}{k} \! \right ) 
\end{equation}
uniformly on compact sets.  Since $\cG$ is $\mathfrak{S}_2$-class, the function $\det^{(2)}(\cG)$ is defined and entire, so the right-hand side of (\ref{eq4}) can be written
\[
 \det^{(2)}(I - z\cG) \prod_j \exp \left ( \frac{\lambda_j(\cG)^2 z^2}{2} + \frac{\lambda_j(\cG)^3 z^3}{3} \right )
\]
Similarly, the left-hand side of (\ref{eq4}) can be written
\[
\psi_n(z) \prod_j \exp \left ( \frac{\lambda_j(A'_n)^2 z^2}{2} + \frac{\lambda_j(A'_n)^3 z^3}{3} \right )
\]
because $\sum_j \lambda_j(A'_n) = \tr(A'_n) = 0$. Since $t(C_3,G_n) \rightarrow t(C_3,\mathcal{G})$, the cubic terms can be cancelled in (\ref{eq4}) and the quadratic terms behave predictably, because $\sum_j \lambda_j(A)^2 = \|\cG\|_2^2$ (the Hilbert-Schmidt norm) and $\|G_n\|_2^2 = 2|E(G_n)|/n^2$:
\[
\psi_n(z) \exp \left (\frac{ z^2 |E(G_n)|}{n^2} \right ) \rightarrow \det^{(2)}(I - z\cG)  \exp \left ( \frac{z^2 \|\cG\|_2^2}{2} \right )
\]
Since the edge density $\sum_j \lambda_j(A')^2 = \|G_n\|_2^2 = 2|E(G_n)|/n^2$ converges to $\|\cG\|_1$, this can be rewritten as
\[
\lim_{n \rightarrow \infty} \psi_n(z) = \det^{(2)}(I - z\cG)  \exp \left ( z^2 \cdot \frac{\|\cG\|_2^2 - \|\cG\|_1}{2} \right ) = \psi_\cG(z)
\]
where the limit can be interpreted as uniform convergence on compact subsets of $\mathbb{C}$ or pointwise.

Then (1) and (2) follow, and, since (2) holds, Cauchy's integral formula implies that (3) holds as well (by integrating $\psi_n(z)/(2 \pi i z^k)$ along a circle around the origin for each $k \geq 0$).  Since $\psi_n(z)$ has only real roots (being the characteristic polynomial of a real symmetric matrix), the limit is of Laguerre-P\'olya class.  That the roots of $\psi_\cG(z)$ are the reciprocals of the eigenvalues of the kernel operator corresponding to $\cG$ now follows from Definitions 1 and 2.
\end{proof}

Theorem \ref{thm:convergence} has immediate consequences from various properties of determinants, for example the following result.  Here we introduce the notation $\cG \oplus_p \cH$ for the {\em $p$-disjoint union} of $\cG$ and $\cH$, the graphon whose kernel $W$ is given by
$$
W(x,y) = \left \{ \begin{array}{ll}
\cG \left ( \frac{x}{p} , \frac{y}{p} \right ) & \text{if } x,y \in [0,p) \\
\cH \left ( \frac{x-p}{1-p}, \frac{y-p}{1-p} \right ) & \text{if } x,y \in [p,1] \\
0 & \text{otherwise} \\
\end{array}
\right . .
$$

\begin{corollary}
    Given two graphons $\cG$ and $\cH$, the graphon $\cG \oplus_p \cH$ which is their disjoint union has the property that
    $$
    \psi_{\cG \oplus_p \cH}(z) = \psi_\cG(pz) \psi_\cH((1-p)z) .
    $$
\end{corollary}
\begin{proof}
    The kernel of $\cG \oplus_p \cH$ is $\cG' + \cH'$, where $\cG'(x,y) = \cG(x/p,y/p)$ and $\cH'(x,y) = \cH((x-p)/(1-p),(y-p)/(1-p))$ (interpreting functions to be zero outside $[0,1] \times [0,1]$).  By \cite{GoKr69} (Section VI.2) and the fact that $\cG' \cH' = 0$ (multiplication interpreted as composition), the Hilbert-Carleman determinant satisfies 
    \begin{align*}
    \det^{(2)}(I - z \cG \oplus \cH) &= \det^{(2)}((I - z \cG')(I- z \cH'))  \\
    &= \det^{(2)}(I - z p\cG) \det^{(2)}(I - z (1-p)\cH) e^{-z \sum_i \lambda_i(\cG' \cH')}.
    \end{align*}
    But $\cG' \cH' = 0$, so $\sum_i \lambda_i(\cG' \cH') = 0$.  Then
    \begin{align*}
        \psi_{\cG \oplus \cH}(z) &= \det^{(2)}(I - z \cG \oplus \cH) \exp \left (z^2 \cdot \frac{\|\cG \oplus \cH \|_2^2 - \|\cG \oplus \cH \|_1^1}{2} \right ) \\
        &= \det^{(2)}(I - z \cG') \det^{(2)}(I - z \cH') \\
        & \qquad \cdot \exp \left [ \frac{z^2}{2} \left ( p^2 \|\cG\|_2^2 + (1-p)^2 \| \cH \|_2^2 - p^2 \|\cG\|_1 - (1-p)^2 \| \cH \|_1 \right ) \right ] \\
        &= \psi_{\cG'}(z) \psi_{\cH'}(z) = \psi_{\cG}(pz) \psi_{\cH}((1-p)z).
    \end{align*}
\end{proof}

One can also obtain the above result by viewing $\cG' \oplus_p \cH'$ as an operator on a direct sum of Hilbert spaces, whose spectrum is easy to describe.  

\section{Special Cases}

Recall that $\tr \cG = \sum_j \lambda_j$, the sum of the eigenvalues of (the kernel of) $\cG$, and that $\cG$ is ``trace class'' (aka ``nuclear'') if this sum converges absolutely. When $\cG$ is trace class, we may write
\begin{equation} \label{eq3}
\psi_\cG(z) = \exp \left(\frac{(\|\cG\|_2^2 - \|\cG\|_1^1) z^2}{2} + z \tr \cG \right) \prod_{s \in \mathcal{S}} \left (1-sz \right ),
\end{equation}
a factorization of the characteristic power series into a monic polynomial-like product whose roots are $\{1/s : s \in \mathcal{S}\}$ where $\mathcal{S}$ is the set of nonzero eigenvalues of $\cG$, and an exponential term.  The quantity $\tr \cG$ is zero if $\cG$ is bipartite: in particular, the eigenvalues comes in pairs $\pm \lambda_j$.  (For more on the spectra of bipartite graphs, see \cite{DoHl19}, in particular Theorem 8.)  In this case, $\psi_G(z)$ has no monomials of odd degree:
$$
\psi_\cG(z) = \exp \left[ (\|\cG\|_2^2 - \|\cG\|_1^1 ) z^2/2 \right ] \prod_{s \in \mathcal{S} \cap \mathbb{R}^+} \left (1-s^2z^2 \right ),
$$
Furthermore, $\|\cG\|_2^2 = \|\cG\|_1^1$ iff $\cG$ is a $0$-$1$ function except for a set of measure zero, as with a simple blow-up of a graph (sometimes called a ``pixel diagram''), so the quadratic term vanishes in the exponential, resulting in
$$
\psi_\cG(z) = \exp \left(z \tr \cG \right) \prod_{s \in \mathcal{S}} \left (1-sz \right ).
$$
We can also use (\ref{eq3}) to obtain a simple expression for the characteristic power series of $p$-quasi-random graphons.
\begin{prop}
A sequence of graphs $G_n$ is $p$-quasirandom iff it converges to a graphon $\cG$ with
\begin{equation} \label{eq2}
\psi_\cG(z) = (1-pz) \exp \left (pz - \frac{p(1-p)}{2} z^2 \right ).
\end{equation}
Furthermore, when this occurs, we can write
\begin{equation} \label{eq7}
\psi_\cG(z) = \sum_{k=0}^\infty z^{k} \sum_{\lambda \in \Lambda(k;i,j)} \frac{(-1)^{j} p^{k-i}}{\eta(\lambda)}
\end{equation}
where $\Lambda(k;i,j)$ is the set of integer partitions of $k$ into $j$ parts of size at least $2$, of which $i$ are of size exactly $2$.
\end{prop}
\begin{proof} 
It is clear from Theorem \ref{thm:QRG} that, if $G_n$ is a $p$-quasirandom sequence, then it converges to a graphon $\cG$ which has only one nonzero eigenvalue $\lambda = p$, and furthermore $(\|\cG\|_2^2 - \|\cG\|_1)/2 = -p(1-p)/2$, so
\begin{equation}\label{eq:QR}
\psi_\cG(z) = (1-pz) \exp \left (pz - \frac{p(1-p)}{2} z^2 \right )
\end{equation}
On the other hand, if (\ref{eq:QR}) holds, then $\cG$ has only one non-zero eigenvalue, and $\psi_\cG$ has only one root, at $z = 1/p$. By the Spectral Theorem, this also holds if and only if the graphon $\cG$, regarded as a function of $[0,1]^2$ up to modification on a set of measure zero, can be written as $p\phi(x)\phi(y)$ for some function $\phi$ normalized so that $\|\phi\|_2=1$. Since the quadratic term in expression (\ref{eq:QR}) is $-p(1-p)/2 = (\|\cG\|_2^2 - \|\cG\|_1)/2$, and $\|\cG\|_2^2 = p^2$, this implies that $\|\cG\|_1 = p$.  But then 
$$
p = \|\cG\|_1 = \|p \phi(x) \phi(y) \|_1 = p \|\phi\|_1^2 \leq p \|\phi\|_2^2 = p
$$
so that $\|\phi\|_1 = \|\phi\|_2$, $\phi$ is the constant function, and $\cG$ is $p$-quasirandom.

To see that $\psi_\cG$ can also be written as in (\ref{eq7}), observe that, by Theorem \ref{thm:QRG} and \ref{thm:convergence}, as well as the discussion following Theorem \ref{thm:HS},
\begin{align*} 
\psi_\cG(z) &= \lim_{n \rightarrow \infty} \sum_{k=0}^n \sum_{H \in \cH_k} z^{k} (-1)^{c(H)} \frac{(\# H \subseteq G_n)}{n^k \eta(H)} \\
&= \sum_{k=0}^\infty \sum_{H \in \cH_k} z^{k} (-1)^{c(H)} \frac{p^{|E(H)|}}{\eta(H)} \\
 &= \sum_{k=0}^\infty z^{k} \sum_{\lambda \in \Lambda(k;i,j)} \frac{(-1)^{j} p^{k-i}}{\eta(\lambda)}
\end{align*}
because partitions of $k$ into parts of size at least $2$ correspond bijectively to elements of $\cH_k$.
\end{proof}

It is also straightforward to see directly that the middle and right-hand expressions in \ref{eq4} are equal.  If $\alpha(z) = \exp(\beta(z))$ is a power series, where the $z^n$ coefficient of $\alpha$ is $a_k$ and the $z^k$ coefficient of $\beta$ is $b_k$, then (by standard facts about exponential generating functions, see, e.g., \cite{Bo17}), letting $\lambda$ be an integer partition with $m_i$ parts of distinct sizes $c_i$, $i = 1$ to $t$,
\begin{align*}
a_k &= \frac{1}{k!} \sum_{\pi \in \Pi} \prod_{B \in \pi} b_{|B|} |B|! \\
&= \frac{1}{k!}  \sum_{\lambda \vdash k} \frac{k!}{\eta(\lambda) \prod_{i=1}^t (c_i-1)!^{m_i}} \prod_{i=1}^t b_{i}^{m_i} c_i!^{m_i} \\
&= \sum_{\lambda \vdash k} \frac{ \prod_{i=1}^t b_{i}^{m_i} c_i^{m_i}}{\eta(\lambda)}
\end{align*}
where $\Pi$ is the set of (set) partitions of an $k$-set.  Thus, if 
\begin{align*}
\beta(z) &= pz - z^2p(1-p)/2 + \log(1-pz) \\
&= pz + \frac{z^2(p^2-p)}{2} - pz - \frac{(pz)^2}{2} - \frac{(pz)^3}{3} - \cdots \\
&= -z^2 \cdot \frac{p}{2} - z^3 \cdot \frac{p^3}{3}  -  z^4 \cdot\frac{p^4}{4} \cdots
\end{align*}
and $\alpha(z) = \exp(\beta)$, then
\begin{align*}
a_k &= \sum_{\lambda \vdash k} \frac{ \prod_{i=1}^t b_{i}^{m_i} c_i^{m_i}}{\eta(\lambda)} \\
&= \sum_{\lambda \in \Lambda(k;r,s)} p^{-r} \frac{ \prod_{i=1}^t (-1/c_i)^{m_i} c_i^{m_i} p^{m_i}}{\eta(\lambda)} \\
&= \sum_{\lambda \in \Lambda(k;r,s)} p^{-r} \frac{ \prod_{i=1}^t (-p)^{m_i}}{\eta(\lambda)} \\
&= \sum_{\lambda \in \Lambda(k;r,s)} (-1)^s p^{k-r} \frac{1}{\eta(\lambda)}.
\end{align*}

\section{Questions}

It is tempting to define instead a characteristic power series without the regularization, i.e., $\det(I - z \cG) = \prod_{j} (1 - \lambda_j z)$, the ``Fredholm determinant''.  However, this product may not converge absolutely.  Using Fourier series, it is straightforward to show that, for the graphon $\cG(x,y)$ defined by 
$$
\cG(x,y) = \left \{ 
\begin{array}{ll}
1 & \text{if $x+y \leq 1$} \\
0 & \text{otherwise}
\end{array}
\right .
$$
we have
$$
\cG(x,y) = \sum_{n = 0}^\infty \frac{4(-1)^n}{(2n+1) \pi} \cos \left (\frac{(2n+1)\pi x}{2} \right ) \cos \left (\frac{(2n+1)\pi y}{2} \right ) 
$$
Since $\{f_n\}_{n=0}^\infty$ with $f_n(t) = \sqrt{2} \cos[(2n+1) \pi t/2]\}$ is an orthonormal family of functions on $[0,1]$, this shows that the spectrum of $\cG$ is $\left \{\frac{2 (-1)^n}{(2n+1)\pi } \right \}_{n=0}^\infty$.  Since the (odd) harmonic series diverges, it follows that $\cG$ is not trace-class, i.e., $\sum_i \lambda_i$ does not converge absolutely.  However, it still converges conditionally, so $\det(I - z \cG)$ is well-defined.

One can also construct a function from $[0,1] \times [0,1] \rightarrow \mathbb{R}^2$ by setting $\cG(x,y) = \sum_{n=0}^\infty \frac{1}{n} (-1)^{\epsilon_n(x) + \epsilon_n(y)}$ where $\epsilon_n(x)$ is the $n$-th digit of $x$ in binary.  Since $\{(-1)^{\epsilon_n}\}_{n=0}^\infty$ is the orthonormal Rademacher system, the spectrum of $\cG$ is the harmonic series.  Thus $\det(I - z \cG)$ is undefined, but because random harmonic series have full support on the real line, the function $\cG$ is not a graphon (and cannot be linearly scaled to become one). 

Thus, we are led to the following question.

\begin{question}
    Does there exist a graphon $\cG$ for which the Fredholm determinant $\det(I - z\cG)$ does not exist?
\end{question}

The characteristic polynomial of the $p$-quasirandom graphon has some possible connections with its probabilistic interpretations, as follows.

\begin{question}
    Let $\cG$ be the $p$-quasirandom graphon.  The function
    $$
    \psi_\cG(z) = (1-pz) \exp(pz - z^2p(1-p)/2)
    $$
    has some unexplained connections with Gaussian probability distributions.  If $M(z)$ is the moment generating function of a normal distribution of mean $p$ and variance $p(1-p)$ -- the normalized limit of a $0$-$1$ random walk with bias $p$ -- then $\psi_\cG(z) = (1-pz)/M(-z)$. Why?
\end{question}

Our next question concerns to what extent some of the above approach can be applied to the many other well-known graph polynomials: matching polynomial, Laplacian characteristic polynomial, chromatic polynomial, etc.  When is it the case that, given some notion of graph convergence, such as left-convergence leading to graphons as above or Benjamini-Schramm convergence of very sparse graphs, these polynomials when suitably normalized converge to some power series?  One motivation for asking this is a related, growing body of work on limits of measures supported on the roots of natural graph polynomials.  For example, building off of work by Sokal \cite{So01} and Borgs-Chayes-Kahn-Lov\'asz \cite{BoChKaLo13}, Ab\'ert-Hubai \cite{AbHu15} showed the convergence of harmonic moments (quantities $\int_K f \, d\nu$ for holomorphic functions $f$ and certain regions $K \subset \mathbb{C}$) of the uniform probability distribution over the chromatic roots of Benjamini-Schramm-convergence graph sequences; subsequently, Csikv\'ari-Frenkel \cite{CsFr16} generalized this to a wide class of graph polynomials (including the characteristic polynomial) and Csikv\'ari-Frenkel-Hladk\'y-Hubai \cite{CsFrHlHu17} showed that it holds even for dense graph (i.e., graphon) convergence with suitable normalization.  From a different perspective, \cite{BuEhFlStTi10} empirically showed that chromatic roots of Erd\H{o}s-R\'enyi random graphs appear to have a scaling limit.

\begin{question}
For which other graph polynomials and graph limit process can the above type of analysis be carried out?  How about for hypergraphs?
\end{question}

Finally, we mention a question that arose in the context of experimentally computing the coefficients of the characteristic power series of that simplest of graphons, the uniform quasirandom graphon.

\begin{question}
It is straightforward to show (by, for example, applying Tur\'an's Inequalities; see \cite{ScPo14}) that the coefficients $c_k$ of $\psi_\cG$ are log-concave for any graphon $\cG$, but the consequences of this for unimodality are unclear because we do not know the sign pattern of the coefficients of $\psi_{\cG}(z)$.  For example, the characteristic power series of a $p=1/2$ quasi-random graphon has sign pattern
$$
+,0,-,+,+,-,-,+,+,+,-,-,+,+,-,-,\ldots
$$
More specifically, can $c_k = 0$ if $k \neq 1$?
\end{question}

\section{Acknowledgments}
Thanks to G.~Clark, C.~Edgar, G.~Fickes, V.~Nikiforov, and A.~Riasanovsky for helpful discussions.  Thanks also to the anonymous referees for their patience and very useful suggestions, especially the identification of a need to use higher-order determinants.


\bibliography{refs}

\begin{thebibliography}{18}
\providecommand{\natexlab}[1]{#1}
\providecommand{\url}[1]{\texttt{#1}}
\expandafter\ifx\csname urlstyle\endcsname\relax
  \providecommand{\doi}[1]{doi: #1}\else
  \providecommand{\doi}{doi: \begingroup \urlstyle{rm}\Url}\fi

\bibitem[Ab\'{e}rt and Hubai(2015)]{AbHu15}
M.~Ab\'{e}rt and T.~Hubai.
\newblock Benjamini-{S}chramm convergence and the distribution of chromatic
  roots for sparse graphs.
\newblock \emph{Combinatorica}, 35\penalty0 (2):\penalty0 127--151, 2015.

\bibitem[Bakan et~al.(2001)Bakan, Craven, and Csordas]{BaCrCs01}
A.~Bakan, T.~Craven, and G.~Csordas.
\newblock Interpolation and the {L}aguerre-{P}\'{o}lya class.
\newblock \emph{Southwest J. Pure Appl. Math.}, \penalty0 (1):\penalty0 38--53,
  2001.

\bibitem[B\'{o}na(2017)]{Bo17}
M.~B\'{o}na.
\newblock \emph{A walk through combinatorics}.
\newblock World Scientific Publishing Co. Pte. Ltd., Hackensack, NJ, 2017.
\newblock An introduction to enumeration and graph theory, Fourth edition [of
  MR1936456], With a foreword by Richard Stanley.

\bibitem[Borgs et~al.(2013)Borgs, Chayes, Kahn, and Lov\'{a}sz]{BoChKaLo13}
C.~Borgs, J.~Chayes, J.~Kahn, and L.~Lov\'{a}sz.
\newblock Left and right convergence of graphs with bounded degree.
\newblock \emph{Random Structures Algorithms}, 42\penalty0 (1):\penalty0 1--28,
  2013.

\bibitem[Chung et~al.(1989)Chung, Graham, and Wilson]{ChuGraWil89}
F.~R.~K. Chung, R.~L. Graham, and R.~M. Wilson.
\newblock Quasi-random graphs.
\newblock \emph{Combinatorica}, 9\penalty0 (4):\penalty0 345--362, 1989.

\bibitem[Csikv\'{a}ri and Frenkel(2016)]{CsFr16}
P.~Csikv\'{a}ri and P.~E. Frenkel.
\newblock Benjamini-{S}chramm continuity of root moments of graph polynomials.
\newblock \emph{European J. Combin.}, 52\penalty0 (part B):\penalty0 302--320,
  2016.

\bibitem[Csikv\'{a}ri et~al.(2017)Csikv\'{a}ri, Frenkel, Hladk\'{y}, and
  Hubai]{CsFrHlHu17}
P.~Csikv\'{a}ri, P.~E. Frenkel, J.~Hladk\'{y}, and T.~Hubai.
\newblock Chromatic roots and limits of dense graphs.
\newblock \emph{Discrete Math.}, 340\penalty0 (5):\penalty0 1129--1135, 2017.

\bibitem[Dole\v{z}al and Hladk\'{y}(2019)]{DoHl19}
M.~Dole\v{z}al and J.~Hladk\'{y}.
\newblock Matching polytons.
\newblock \emph{Electron. J. Combin.}, 26\penalty0 (4):\penalty0 Paper No.
  4.38, 33, 2019.

\bibitem[Gohberg and Kre\u{\i}n(1969)]{GoKr69}
I.~C. Gohberg and M.~G. Kre\u{\i}n.
\newblock \emph{Introduction to the theory of linear nonselfadjoint operators}.
\newblock Translated from the Russian by A. Feinstein. Translations of
  Mathematical Monographs, Vol. 18. American Mathematical Society, Providence,
  R.I., 1969.

\bibitem[Griffin et~al.(2019)Griffin, Ono, Rolen, and Zagier]{GrOnRoZa19}
M.~Griffin, K.~Ono, L.~Rolen, and D.~Zagier.
\newblock Jensen polynomials for the {R}iemann zeta function and other
  sequences.
\newblock \emph{Proc. Natl. Acad. Sci. USA}, 116\penalty0 (23):\penalty0
  11103--11110, 2019.

\bibitem[Harary(1962)]{Ha62}
F.~Harary.
\newblock The determinant of the adjacency matrix of a graph.
\newblock \emph{SIAM Rev.}, 4:\penalty0 202--210, 1962.

\bibitem[Hirschman and Widder(1955)]{HiWi55}
I.~I. Hirschman and D.~V. Widder.
\newblock \emph{The convolution transform}.
\newblock Princeton University Press, Princeton, N. J., 1955.

\bibitem[Lov\'{a}sz(2012)]{Lo12}
L.~Lov\'{a}sz.
\newblock \emph{Large networks and graph limits}, volume~60 of \emph{American
  Mathematical Society Colloquium Publications}.
\newblock American Mathematical Society, Providence, RI, 2012.

\bibitem[Nikiforov(2007)]{Ni07}
V.~Nikiforov.
\newblock Graphs and matrices with maximal energy.
\newblock \emph{J. Math. Anal. Appl.}, 327\penalty0 (1):\penalty0 735--738,
  2007.

\bibitem[Schur and P\'{o}lya(1914)]{ScPo14}
J.~Schur and G.~P\'{o}lya.
\newblock \"{U}ber zwei {A}rten von {F}aktorenfolgen in der {T}heorie der
  algebraischen {G}leichungen.
\newblock \emph{J. Reine Angew. Math.}, 144:\penalty0 89--113, 1914.

\bibitem[Sokal(2001)]{So01}
A.~D. Sokal.
\newblock Bounds on the complex zeros of (di)chromatic polynomials and
  {P}otts-model partition functions.
\newblock \emph{Combin. Probab. Comput.}, 10\penalty0 (1):\penalty0 41--77,
  2001.

\bibitem[Szegedy(2011)]{Sze11}
B.~Szegedy.
\newblock Limits of kernel operators and the spectral regularity lemma.
\newblock \emph{European J. Combin.}, 32\penalty0 (7):\penalty0 1156--1167,
  2011.

\bibitem[Van~Bussel et~al.(2010)Van~Bussel, Ehrlich, Fliegner, Stolzenberg, and
  Timme]{BuEhFlStTi10}
F.~Van~Bussel, C.~Ehrlich, D.~Fliegner, S.~Stolzenberg, and M.~Timme.
\newblock Chromatic polynomials of random graphs.
\newblock \emph{J. Phys. A}, 43\penalty0 (17):\penalty0 175002, 12, 2010.

\end{thebibliography}

\end{document}